\documentclass[11pt]{amsart}
\usepackage{amssymb, amsthm, epsfig}
\usepackage[mathscr]{euscript}

\newtheorem{thm}{Theorem}[section]
\newtheorem{cor}[thm]{Corollary}
\newtheorem{lem}[thm]{Lemma}
\newtheorem{prop}[thm]{Proposition}

\newtheorem{defnn}[thm]{Definition}
\newenvironment{definition}{\begin{defnn} \em}{\end{defnn}}
\newtheorem{remarkk}[thm]{Remark}
\newenvironment{remark}{\begin{remarkk} \em}{\end{remarkk}}
\newtheorem{examplee}[thm]{Example}

\newcommand{\scra}{\mathcal{A}}

\newcommand{\scrc}{\mathcal{C}}
\newcommand{\scrf}{\mathcal{F}}
\newcommand{\scrk}{\mathcal{K}}
\newcommand{\scrr}{\mathcal{R}}
\newcommand{\bbz}{\mathbb{Z}}
\newcommand{\bbr}{\mathbb{R}}
\newcommand{\bbq}{\mathbb{Q}}
\newcommand{\bbn}{\mathbb{N}}

\newcommand{\gu}{\Gamma^U_n}
\newcommand{\imra}{\looparrowright}
\newcommand{\Hom}{\operatorname{Hom}}
\newcommand{\Ker}{\operatorname{Ker}}

\title
[Filtration of the classical knot concordance group]
{Filtration of the classical knot concordance group and  Casson-Gordon invariants}
\author{Taehee Kim}
\date{\today}
\begin{document}
\begin{abstract} It is known that if any prime power branched
cyclic cover of a knot in $S^3$ is a homology sphere, then the
knot has vanishing Casson-Gordon invariants. We construct
infinitely many examples of (topologically) non-slice knots in
$S^3$ whose prime power branched cyclic covers are homology
spheres.  We show that these knots generate an infinite rank
subgroup of $\scrf_{(1.0)}/\scrf_{(1.5)}$ for which Casson-Gordon
invariants vanish in Cochran-Orr-Teichner's filtration of the
classical knot concordance group . As a corollary, it follows that
Casson-Gordon invariants are not a complete set of obstructions to
a second layer of Whitney disks.
\end{abstract}

\maketitle \vspace*{-2em}

\section{Introduction}
\label{sec:introduction} A knot in the 3-sphere is (topologically)
slice if it bounds a locally flat 2-disk in the 4-ball. Two knots
are said to be (topologically) concordant if the connected sum of
one and the mirror image of the other with reversed orientation is
slice. (Equivalently, there is a locally flat embedding of an
annulus $S^1\times [0,1]$ into $S^3\times [0,1]$ whose
restrictions to the boundary components give the knots.) This
concordance relation is an equivalence relation, and the
concordance classes form an abelian group $\scrc$, the classical
knot concordance group, under the connected sum operation. In
$\scrc$, the identity element is the class of slice knots.

In \cite{COT1}, Cochran, Orr, and Teichner (henceforth COT) define
a geometric filtration of the classical knot concordance group
$\scrc$
\[
 0\subset\cdots\subset\scrf_{(n.5)}\subset\scrf_{(n)}\subset\cdots
 \subset\scrf_{(1.5)}\subset\scrf_{(1.0)}\subset\scrf_{(0.5)}
 \subset\scrf_{(0)}\subset\scrc
\]
where $\scrf_{(m)}$ is the set of $(m)$-solvable knots. (See
Definition \ref{def:(n)-solvability}.) They show that
$(1.5)$-solvable knots have vanishing Casson-Gordon invariants and
that $\scrf_{(2.0)}/\scrf_{(2.5)}\neq 0$, thus giving the first
examples of knots with vanishing Casson-Gordon invariants which
are not (topologically) slice. (Refer to \cite{CG} for
Casson-Gordon invariants.) In \cite{COT2}, they extend their
results to show $\scrf_{(2.0)}/\scrf_{(2.5)}$ has infinite rank.
We improve their results further and prove:

\begin{thm}[Main Theorem]
\label{thm:main theorem} In the above filtration,
$\scrf_{(1.0)}/\scrf_{(1.5)}$ has an infinite rank subgroup of
knots for which Casson-Gordon invariants vanish.
\end{thm}

\noindent Theorem~\ref{thm:main theorem} implies that
Casson-Gordon invariants are not a complete set of obstructions to
$(1.5)$-solvability. By contrast to the above result, the examples
of \cite{COT1} are $(2.0)$-solvable.

To show Casson-Gordon invariants of our examples vanish, we use
the following theorem of Livingston.

\begin{thm}
\label{thm:homology sphere} (\cite[Theorem 0.5]{Liv}) A knot $K$
has a prime power branched cyclic cover with nontrivial homology
if and only if its Alexander polynomial has a nontrivial factor
that is not an n-cyclotomic polynomial with n divisible by three
distinct primes.
\end{thm}

\noindent The group of examples in Theorem~\ref{thm:main theorem}
have a spanning set of knots with a fixed Seifert form and
Alexander polynomial. The shared Alexander polynomial of these
generators is $(\Phi_{30})^2$, the square of the 30-cyclotomic
polynomial. By Theorem~\ref{thm:homology sphere}, these generators
have {\it prime power} branched cyclic covers which are homology
spheres. (One might compare this to the fact that if every finite
branched cyclic cover of a knot is a homology sphere, then its
Alexander polynomial is 1, hence the knot is topologically slice
by Freedman's work \cite{F}.) It follows from the definition of
Casson-Gordon invariants that Casson-Gordon invariants vanish for
these knots. (See Proposition~\ref{prop:vanishing CG}.) In fact,
since any prime power branched cyclic cover is a homology sphere
all the concordance invariants known prior to
Cochran-Orr-Teichner's $L^{(2)}$-signature invariants, such as
Gilmer's extension of Casson-Gordon invariants (\cite{G}), Kirk
and Livingston's twisted Alexander invariants (\cite{KL}), and
Letsche's invariants (\cite{Le}), vanish for these knots.

Theorem \ref{thm:main theorem} has a significant geometric
consequence. Freedman's disk embedding theorem (\cite{F}),
together with the Cappell-Shaneson homology surgery approach
(\cite{CS}) to classifying knot concordance group, suggest that
the Casson-Gordon invariants obstruct the construction of a second
layer of Whitney disks for a Cappell-Shaneson surgery kernel of an
algebraically slice knot. That this is so was shown in
\cite[Section 8 and 9]{COT1}. Indeed, in \cite{COT1}, they showed
that a knot is $(1.5)$-solvable if and only if for zero surgery on
the knot in $S^3$, there exists an $H_1$-bordism which contains a
spherical Lagrangian admitting a Whitney tower of height $(1.5)$.
(See \cite[Theorem 8.4]{COT1} and Section \ref{sec:filtration} in
this paper). Since $(1.5)$-solvable knots have vanishing
Casson-Gordon invariants, it follows that Casson-Gordon invariants
obstruct a Whitney tower of height $(1.5)$ in the above sense.
Precise definitions of a Whitney tower and other terminologies are
given in \cite[Section 8]{COT1} and are reviewed in Section
\ref{sec:filtration} in this paper.

We briefly discuss Whitney towers here. In 4-manifolds, Whitney
disks may no longer be embedded, but may themselves have
intersections, which might or might not occur in algebraically
cancelling pairs. If these intersections occur in algebraically
cancelling pairs, one can construct immersed Whitney disks for
these cancelling pairs of points in the usual manner. Very roughly
speaking, a Whitney tower is obtained by iterating this procedure.
We have the following corollary of Theorem \ref{thm:main theorem}.

\begin{cor}
\label{cor:1.5 and Whitney towers} There is an algebraically slice
knot with vanishing Casson-Gordon invariants such that zero
surgery on the knot in $S^3$ does not bound an $H_1$-bordism which
contains a spherical Lagrangian admitting a Whitney tower of
height (1.5).
\end{cor}
\begin{proof}
It follows from Theorem~\ref{thm:main theorem} and \cite[Theorem
8.4]{COT1}.
\end{proof}

\noindent Corollary \ref{cor:1.5 and Whitney towers} says that
Casson-Gordon invariants are not a complete set of obstructions to
a second layer of Whitney disks.

To find the knots generating the subgroup in Theorem~\ref{thm:main
theorem}, we follow the method of COT. We begin by constructing a
ribbon knot with the rational Alexander module
$\bbq[t,t^{-1}]/(\Phi_{30}(t))^2$. In particular, its Alexander
polynomial is $(\Phi_{30}(t))^2$. Henceforth we refer to this
ribbon knot as the {\it seed knot} to the examples of Theorem
\ref{thm:main theorem}. (See Remark~\ref{rem:seed knot}. See
\cite{K} for the definition of a ribbon knot. In particular, a
ribbon knot is a slice knot.) We modify this seed knot using a
family of Arf invariant zero knots in a way described in
\cite[Setion 3]{COT2} and reviewed in Section~\ref{sec:n-solvable
knots} in this paper. The resulting knots are shown to have the
same Seifert form with the seed knot, so their prime power
branched cyclic covers are also homology spheres by
Theorem~\ref{thm:homology sphere}. Another important fact, which
will be used significantly in this paper, is that
$\bbq[t,t^{-1}]/(\Phi_{30}(t))^2$ has a unique nontrivial proper
submodule. (See the proofs of Lemma~\ref{lem:uniqueness},
Proposition~\ref{prop:1-solvability} and Theorem~\ref{thm:main
theorem}.)

This paper is organized in the following manner. In
Section~\ref{sec:seed knot}, we construct a ribbon knot whose
rational Alexander module is cyclic of order $(\Phi_{30}(t))^2$,
i.e., $\bbq[t,t^{-1}]/(\Phi_{30}(t))^2$. This will be our seed
knot. In Section~\ref{sec:filtration}, we explain the definition
and properties of the Cochran-Orr-Teichner filtration of the
classical knot concordance group and its relation to Whitney
towers. In Section~\ref{sec:n-solvable knots}, we discuss how to
construct a family of $(n)$-solvable knots from a given ribbon
knot using a certain Arf invariant zero knot. This method, applied
to the ribbon knot mentioned above, will be used to construct the
generators of the desired subgroup. In
Section~\ref{sec:detecting}, $L^{(2)}$-signatures and their
properties are reviewed. Finally, in Section~\ref{sec:main
theorem}, we provide the construction of a set of generators of
the subgroup in Theorem~\ref{thm:main theorem} and the proof of
Theorem~\ref{thm:main theorem}.

\begin{remark}
\label{rem:general n} For any $n\in\bbn$ which is divisible by at
least three distinct primes, we can also find an infinite rank
subgroup of $\scrf_{(1.0)}/\scrf_{(1.5)}$ that has generators with
the Alexander polynomial $(\Phi_{n}(t))^2$ (the square of the
$n$-cyclotomic polynomial). The only difference in the proof will
be finding a seed knot with the rational Alexander module
$\bbq[t,t^{-1}]/(\Phi_{n}(t))^2$ as we do for $n = 30$ in
Section~\ref{sec:seed knot}.
\end{remark}

The author thanks Kent Orr for his helpful advice and
encouragement. He also would like to thank Chuck Livingston,
Se-Goo Kim, and Jae Choon Cha for helpful conversations.

\section{construction of the seed knot}
\label{sec:seed knot} In this section, we construct our seed knot.
That is, we will construct a knot which is a ribbon knot and has
the rational Alexander module $\bbq[t,t^{-1}]/(\Phi_{30}(t))^2$.
Then by Theorem~\ref{thm:homology sphere} of C.~Livingston, the
seed knot will have prime power branched cyclic covers that are
homology spheres.

First, we find a Seifert matrix whose Alexander polynomial is
$\Phi_{30}(t)$. Recall that $\Phi_{30}(t) = t^8 + t^7 - t^5 - t^4
- t^3 + t + 1$. This can be done by applying Levine's arguments in
\cite[14 on page 236]{L} which originated from Seifert~\cite{S}.
But we will need a Seifert surface that is a boundary connected
sum of disks with two bands, and it's not clear how to find such a
Seifert surface from the resulting Seifert matrix. So we modify
Levine's arguments a little. The final matrix is
\[
A=\left(\begin{matrix}
 0 & 1 & 1 & 0 & 0 & 0 & 0 & 0 \cr
 0 & 0 & 1 & 0 & 0 & 0 & 0 & 0 \cr
 1 & 1 & 0 & 1 & -9 & 0 & 0 & 0 \cr
 0 & 0 & 0 & 0 & 1 & 0 & 0 & 0 \cr
 0 & 0 & -9 & 1 & 0 & 1 & 26 & 0 \cr
 0 & 0 & 0 & 0 & 0 & 0 & 1 & 0 \cr
 0 & 0 & 0 & 0 & 26 & 1 & 24 & 1 \cr
 0 & 0 & 0 & 0 & 0 & 0 & 0 & 1
\end{matrix}\right)
\]

\begin{figure}[htb]
\includegraphics[scale=.5]{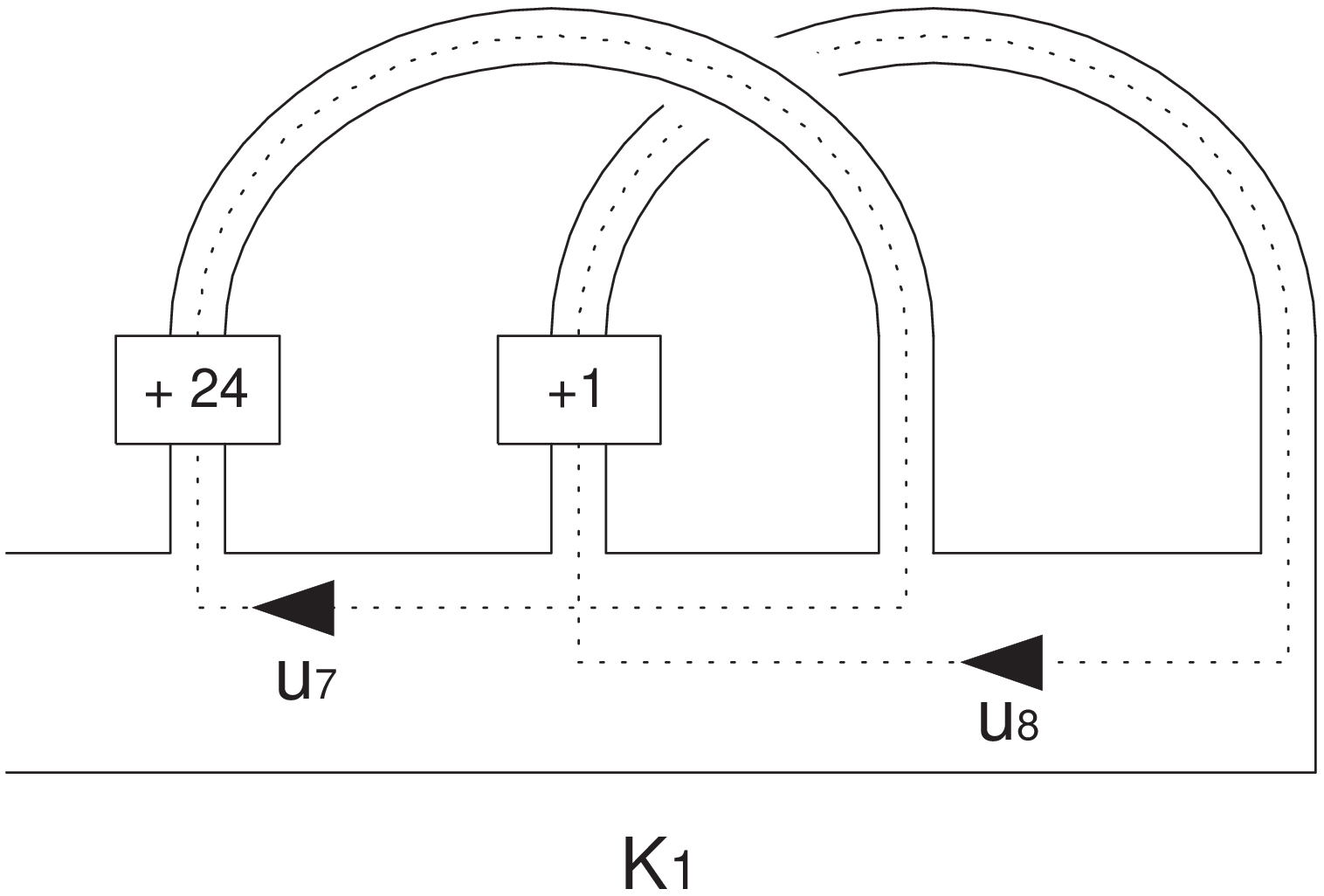}
\caption{}\label{fig:K_1}
\end{figure}

\noindent That is, $\det(A^T-tA) = \Phi_{30}(t)$. One can easily
construct a knot, say $K_1$, and its Seifert surface whose Seifert
matrix with respect to a certain choice of basis is the matrix
$A$. The Seifert surface is obtained in the usual way as the
boundary connected sum of disks with two bands under proper twists
and intertwining among bands. Figure~\ref{fig:K_1} is a part of
$K_1$ and its Seifert surface. The rectangles containing integers
symbolize full twists between the two strands which pass
vertically through the rectangles. Thus the rectangle labelled
$+24$ symbolizes 24 right-handed full twists. Let $u_i, 1\le i\le
8$, be the simple closed curves on the Seifert surface each of
which goes once around a band. With proper orientations,
$\{u_i\}_{1\le i\le 8}$ is a basis with respect to which the
Seifert matrix is the matrix $A$. It is known that $A^T - tA$ is a
presentation matrix of the rational Alexander module of $K_1$. By
column and row operations on $A^T - tA$ over
$\bbq[t,t^{-1}]$-coefficients, we determine that the rational
Alexander module of $K_1$ is isomorphic to
$\bbq[t,t^{-1}]/\Phi_{30}(t)$, whose only generator is represented
by a dual of $u_8$. (A dual of $u_8$ is a simple closed curve in
the complement of the Seifert surface such that it has linking
number one with $u_8$ and no linking with the other $u_i$'s.)

Let $K_2 = K_1\#(-K_1)$, the connected sum of $K_1$ and its
inverse. Then $K_2$ is a ribbon knot. (See, for instance,
\cite[Proposition 5.10 p.83]{K}.) Its rational Alexander module is
$\bbq[t,t^{-1}]/\Phi_{30}(t) \oplus \bbq[t,t^{-1}]/\Phi_{30}(t)$,
and its Seifert surface is obtained as the boundary connected sum
of the Seifert surface of $K_1$ and that of $-K_1$ which is the
mirror image of the Seifert surface of $K_1$. See
Figure~\ref{fig:k2} below. Let $M_2$ denote zero surgery on $K_2$
in $S^3$. The rational Alexander module of $K_2$,
$H_1(M_2;\bbq[t,t^{-1}])$, is generated by $v_i, 1\le i \le 16$,
where for $1\le i\le 8$, $v_i = u_i$, and for $9\le i\le 16$,
$v_i$ is the mirror image of $-u_{(17-i)}$. With this choice of
basis, the Seifert matrix of $K_2$ is the matrix $B =
\left(b_{ij}\right), 1\le i,j\le 16$, defined by
 \[
 b_{ij} = \left\{ \begin{array}
            {r@{\quad:\quad}l}
            a_{ij} & 1\le i,j\le 8 \\
            -a_{(17-i)(17-j)} & 9\le i,j\le 16 \\
            0 & \mbox{otherwise}
            \end{array} \right.
 \]

\begin{figure}[htb]
\includegraphics[scale=0.9]{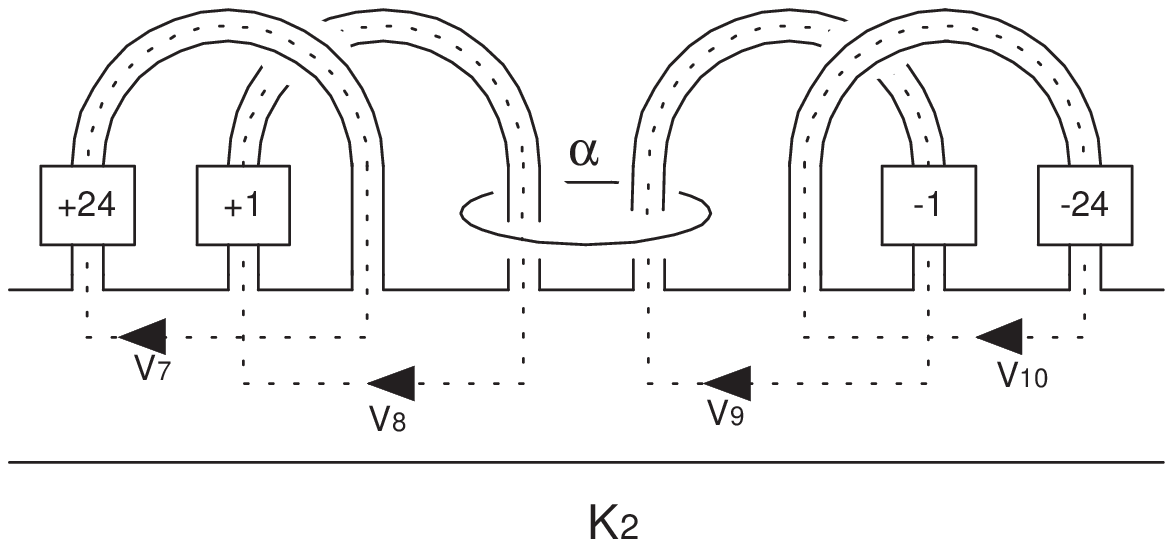}
\caption{}\label{fig:k2}
\end{figure}

\begin{figure}[htb]
\includegraphics[scale=.9]{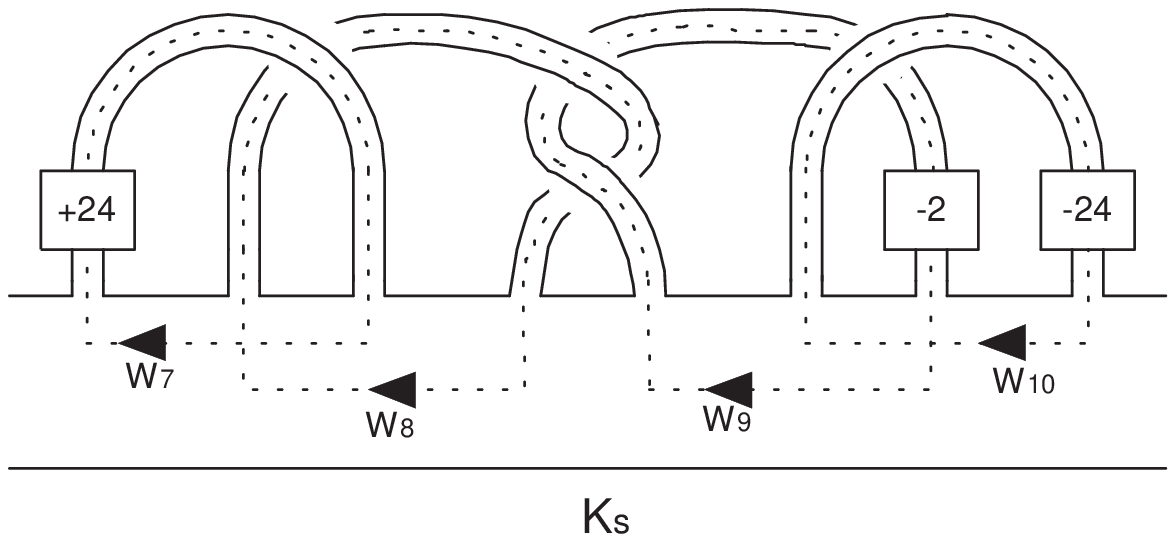}
\caption{}\label{fig:ks}
\end{figure}

Even though $K_2$ is a ribbon knot, its rational Alexander module
is generated by two elements. In particular, it's not cyclic. So
we modify $K_2$ a little more. Choose an unknot $\alpha$ around
the Seifert surface of $K_2$ as in Figure~\ref{fig:k2}. After $+1$
surgery on $\alpha$, $K_2$ will be modified to a new knot, say
$K_s$, in $S^3$ since the resulting ambient manifold obtained by
$+1$ surgery on an unknot in $S^3$ is homeomorphic with $S^3$. A
part of $K_s$ is illustrated in Figure~\ref{fig:ks}. Let $w_i,
1\le i\le 16$, denote the image of $v_i$ under the surgery.
$\{w_i\}_{1\le i\le 16}$ is a basis of the Seifert form of $K_s$.
The Seifert matrix with respect to this basis is obtained by
changing the matrix $B$ such that only $b_{ij}$ with $7\leq i,j
\leq 10$ are changed from

\centerline{ $\left(\begin{matrix} 24 & 1 & 0 & 0\cr
 0 & 1 & 0 & 0\cr
 0 & 0 & -1 & 0 \cr
 0 & 0 & -1 & -24
\end{matrix}\right)$ \hspace{1cm}to\hspace{1cm}
$\left(\begin{matrix} 24 & 1 & 0 & 0\cr
 0 & 0 & 1 & 0\cr
 0 & 1 & -2 & 0 \cr
 0 & 0 & -1 & -24
\end{matrix}\right)$.
}
\noindent Denote the resulting matrix by $C$. Let $M_s$ denote
zero surgery on $K_s$ in $S^3$.

\begin{prop}
\label{prop:Alexander module of K_s} The rational Alexander module
of $K_s$ is cyclic of order $(\Phi_{30}(t))^2$, {\em i.e.},
$H_1(M_s;\bbq[t,t^{-1}])\cong\bbq[t,t^{-1}]/(\Phi_{30}(t))^2$.
\end{prop}
\begin{proof}
$C^T - tC$ is a presentation matrix of $H_1(M_s;\bbq[t,t^{-1}])$.
By column and row operations on $C^T - tC$ over
$\bbq[t,t^{-1}]$-coefficients, we find out that
$H_1(M_s;\bbq[t,t^{-1}])\cong\bbq[t,t^{-1}]/(\Phi_{30}(t))^2$
whose only generator is a dual of $w_9$.
\end{proof}

\begin{prop}
\label{prop:seed knot is ribbon} $K_s$ is a ribbon knot.
\end{prop}
\begin{proof}
One can construct a ribbon disk for the ribbon knot $K_2$ using
the method in \cite[Proposition 5.10]{K}. In particular, a ribbon
disk can be obtained such that $\alpha$ is disjoint from the
ribbon disk and the spanning disk of $\alpha$ has no intersection
with the singularities of the ribbon disk. After $+1$ surgery
along $\alpha$, the image of the ribbon disk of $K_2$ would be a
ribbon disk of $K_s$.
\end{proof}

\begin{remark}
\label{rem:seed knot} By Proposition \ref{prop:Alexander module of
K_s} and \ref{prop:seed knot is ribbon}, $K_s$ has the rational
Alexander module which is cyclic of order $(\Phi_{30}(t))^2$ and
it is a ribbon knot. So $K_s$ is our desired seed knot.
\end{remark}

\section{filtering the knot concordance group and Whitney towers}
\label{sec:filtration} This section and the next two sections are
brief expositions of some of the work in \cite{COT1} and
\cite{COT2}. These sections contain no new results but serve to
clarify ideas and make this paper more self-contained.

In \cite{COT1}, Cochran, Orr, and Teichner established a geometric
filtration of the knot concordance group $\scrc$
\[
 0\subset\cdots\subset\scrf_{(n.5)}\subset\scrf_{(n)}\subset\cdots
 \subset\scrf_{(1.5)}\subset\scrf_{(1.0)}\subset\scrf_{(0.5)}
 \subset\scrf_{(0)}\subset\scrc
\]
where $\scrf_{(n)}$ is the subgroup of \emph{(n)-solvable knots}
for $n\in\{0, 0.5, 1.0, 1.5, \cdots\}$. The precise definition of
this filtration and Whitney towers, and their relations will be
discussed in this section.

Let $G^{(i)}$ denote the {\em i-th derived subgroup} of a group
$G$, inductively defined by $G^{(0)} \equiv G$ and $G^{(i+1)}
\equiv [G^{(i)},G^{(i)}]$. For a CW-complex $W$, we denote the
regular covering of $W$ corresponding to the subgroup
$\pi_1(W)^{(n)}$ by $W^{(n)}$. If $W$ is a spin $4$-manifold, then
we have the usual intersection form
\[
 \lambda_n : H_2(W^{(n)})\times H_2(W^{(n)})\longrightarrow
 \bbz[\pi_1(W)/\pi_1(W)^{(n)}]
\]
A more detailed description of $\lambda_n$ and the
self-intersection invariant $\mu_n$ can be found in \cite[Chapter
5]{W} and \cite[Section 7]{COT1}. In particular, $\lambda_0$ is
the ordinary intersection form on $H_2(W)$.

Now fix a closed oriented 3-manifold $M$.
\begin{definition}
\label{def:H_1-bordism} An $H_1$-{\em bordism} is a 4-dimensional
spin manifold $W$ with boundary $M$ such that the inclusion map
induces an isomorphism $H_1(M)\xrightarrow{\cong}H_1(W)$.
\end{definition}

An {\em $(n)$-surface} is a generic immersion of a closed oriented
surface $F$, say $f : F\imra X$, such that $f_*(\pi_1(F))\le
\pi_1(X)^{(n)}$.

\begin{definition}
\label{def:n-Lagrangian} Let $W$ be an $H_1$-bordism such that
$\lambda_0$ is a hyperbolic form on $H_2(W)$.
\begin{enumerate}
\item A {\em Lagrangian} for $\lambda_0$ is a direct summand of
$H_2(W)$ of half rank on which $\lambda_0$ vanishes.
\item An {\em $(n)$-Lagrangian} is a submodule of $H_2(W^{(n)})$
on which $\lambda_n$ and $ \mu _n$ vanish and which maps onto a
Lagrangian of $\lambda_0$ on $H_2(W)$.
\item   A {\em spherical Lagrangian} is a submodule of
$\pi_2(W)$ on which $ \lambda_n, \mu_n$ $(n\ge 0)$ vanish and
which maps onto a Lagrangian of $\lambda_0$.
\item For $k\le n$, {\em $(k)$-duals} of an $(n)$-Lagrangian
generated by $(n)$-surfaces $\ell_1, \cdots, \ell_g$ are
$(k)$-surfaces $d_1, \cdots, d_g$ such that $H_2(W)$ has rank $2g$
and
\[
\lambda_k(\ell_i,d_j) = \delta_{i,j}.
\]

\end{enumerate}
\end{definition}

Before giving the definition of $(n)$-solvability, we discuss
Whitney towers. Let $W$ be a 4-manifold with boundary $M$ and
$\gamma$ be a framed circle in $M$. A {\it Whitney disk} is an
immersed disk $\Delta$ in $W$ which bounds $\gamma$ and such that
the unique framing on the normal bundle of $\Delta$ restricts to
the given framing on $\gamma$. $\gamma$ is called its {\it Whitney
circle}.

\begin{definition}
\label{def:whitney tower}
\begin{enumerate}
\item A {\em Whitney tower of height $(0)$} is a collection $\scrc_0$ of 2-spheres
$S_i\imra W^4$.
\item For $n\in\bbn$, a {\em Whitney tower of height $(n)$ on
$\scrc_0$} is a sequence $\scrc_j = \{\Delta_{j,k}\}_k, j =
1,\dots,n$, of collections of framed immersed Whitney disks
$\Delta_{j,k}$ in general position such that for $j=2,\dots,n$,
the collection $\scrc_j$ pairs up all
$\scrc_{j-1}$-(self)-intersections and has interiors disjoint from
$\scrc_1,\dots,\scrc_{j-1}$.
\item For $n\in\bbn$, a {\em Whitney tower of height $(n.5)$ on $\scrc_0$}
is a sequence $\scrc_j = \{\Delta_{j,k}\}_k, j = 1,\dots,n+1$ of
collections of framed immersed Whitney disks such that $\scrc_1,
\dots, \scrc_n$ consist of a Whitney tower of height $n$ on
$\scrc_0$ and $\scrc_{n+1}$ pairs up all
$\scrc_{n}$-(self)-intersections and has interiors disjoint from
$\scrc_1,\dots, \scrc_{n-1}$ (but $\scrc_{n+1}$ is allowed to
intersect the previous collection $\scrc_{n}$).
\end{enumerate}
\end{definition}
Refer to \cite[Section 7]{COT1} for more details about Whitney
towers.

\begin{definition}
\label{def:(n)-solvability} A $3$-manifold $M$ is
\emph{$(n)$-solvable} (resp. \emph{$(n.5)$-solvable}) if there is
an $H_1$-bordism $W$ which contains an $(n)$-Lagrangian (resp.
$(n+1)$-Lagrangian) with $(n)$-duals. If $M$ is zero surgery on a
knot or a link then the corresponding knot or link is called
$(n)$-solvable (resp. $(n.5)$-solvable).
\end{definition}

\noindent In Definition \ref{def:(n)-solvability}, $M$ is said to
be {\em (n)-solvable via} (resp. {\em (n.5)-solvable via}) $W$,
and $W$ is called an {\em (n)-solution} (resp. {\em
(n.5)-solution}) for $M$.

\begin{thm}
\label{thm:(n)-solvability} (\cite[Theorem 8.4, 8.8]{COT1}) Let
$M$ be a closed oriented $3$-manifold and $n\in\{0, 0.5, 1.0, 1.5,
\cdots\}$. Then $M$ is $(n)$-solvable if and only if there is an
$H_1$-bordism which contains a spherical Lagrangian admitting a
Whitney tower of height~$(n)$.
\end{thm}

\begin{remark}
\label{rem:n-solvability} The exterior of a slice disk is an
$(n)$-solution for the slice knot (and for its zero surgery $M$)
for all $n$.
\end{remark}

\section{constructing $(n)$-solvable knots}
\label{sec:n-solvable knots} In this section, we obtain an
$(n)$-solvable knot by modifying a given ribbon knot $K$. For this
purpose, we make use of a grafting construction, which produces a
satellite knot of $K$. For further details on this construction
and more general cases, the reader should consult \cite[Section
3]{COT2}.

Simply speaking, seize a collection of parallel strands of $K$ in
one hand and tie these into a knot, say $J$. More precisely,
choose a circle, say $\eta$, in $S^3\setminus K$ which bounds an
embedded disk in $S^3$. Now cut open $K$ along this disk and tie
all the strands passing through this disk into $J$, or more
exactly, through a tubular neighborhood of $J$ with $0$-framing.
Then the resulting ambient manifold is still homeomorphic with
$S^3$, and under this identification, we obtain a new knot $K'$
which is the image of $K$. We denote the resulting knot $K'$ by
$K(J,\eta)$. Moreover, this construction has another very useful
description. $K(J,\eta)$ is obtained by taking the union of the
exterior of $\eta$ and that of $J$ along the boundary in such a
way that the resulting ambient manifold is homeomorphic with
$S^3$.

Making use of the above construction, we get the following
proposition due to COT. We outline a proof here for completeness
and to establish notation for what follows. $M$ (resp. $M_J$)
denotes zero surgery on $K$ (resp. $J$) in $S^3$. Note that a knot
is $(0)$-solvable if and only if it has Arf invariant zero.
(\cite[Remark 8.2]{COT1}.)
\begin{prop}
\label{prop:n-solvability} If $\eta\in \pi_1(M)^{(n)}$ and $J$ has
Arf invariant zero, then $K(J,\eta)$ is $(n)$-solvable.
\end{prop}
\begin{proof}
This is a special case of \cite[Proposition 3.1] {COT2}. Let $W$
be the exterior of a ribbon disk for $K$ in $B^4$. (Note that $W$
may be viewed as an $(n)$-solution.) Let $W_J$ be the
$(0)$-solution for $J$ such that a canonical epimorphism
$\pi_1(M_J)\longrightarrow\bbz$ extends to $\pi_1(W_J)$. By doing
surgery on elements in $\pi_1(W_J)^{(1)}$, we can assume that
$\pi_1(W_J) \cong \bbz$. Let $\mu_J$ denote the meridian of a
tubular neighborhood of $J$ and let $\ell_J$ be the $0$-framed
longitude. Then $\partial W_J = M_J = E_J\cup(S^1\times D^2)$
where $S^1\times \{*\}$ is $\mu_J$, and $\{*\}\times\partial D^2$
is $\ell_J$. Let $W'$ be the 4-manifold obtained from $W_J$ and
$W$ by identifying the solid torus $S^1\times D^2\subset\partial
W_J$ with $\eta\times D^2\subset\partial W$. Observe that
$\partial W' = M'$, zero surgery on $K' = K(J,\eta)$. Then $W'$ is
an $(n)$-solution for $K'$. See \cite[Proposition 3.1]{COT2} for
more details.
\end{proof}

\section{Detecting $(n)$-solvability using $L^{(2)}$-signatures}
\label{sec:detecting} A group $\Gamma$ is called {\em
poly-torsion-free-abelian (PTFA)} if it admits a normal series
$\left<1\right>=G_0\triangleleft
G_1\triangleleft\dots\triangleleft G_n=\Gamma$ such that the
factors $G_{i+1}/G_i$ are torsion-free abelian. If $\Gamma$ is
PTFA, then the group ring $\bbq\Gamma$ is a right Ore domain,
hence $\bbq\Gamma$ embeds in its classical right ring of quotients
$\scrk_\Gamma$. (\cite[Proposition 2.5]{COT1}.) Let $M$ be an
oriented closed 3-manifold. Suppose $\phi :
\pi_1(M)\longrightarrow\Gamma$ is a homomorphism where $\Gamma$ is
a PTFA group and suppose there are an oriented compact
$4$-manifold $W$ bounded by $M$ and a homomorphism $\psi :\pi_1(W)
\longrightarrow\Gamma$ which extends $\phi$, \emph{i.e.},
$(M,\phi) = \partial(W,\psi)$. Then the \emph{(reduced)
$L^{(2)}$-signature} or \emph{von Neumann $\rho$-invariant}
$\rho(M,\phi)\in\bbr$ is defined to be $\rho(M,\phi) =
\sigma^{(2)}_\Gamma(W,\psi) - \sigma_0(W)$ where
$\sigma^{(2)}_\Gamma$ is the $L^{(2)}$-signature of the
intersection form on $H_2(W;\scrk_\Gamma)$ and $\sigma_0$ is the
ordinary signature. We refer the reader to \cite[Section 5]{COT1}
for more discussion of $L^{(2)}$-signatures. The following
theorem, due to COT, gives an obstruction for a knot being
$(n.5)$-solvable.

\begin{thm}
\label{thm:facts of rho invariant}(\cite[Theorem 4.2]{COT1})
Suppose $\Gamma$ is an $(n)$-solvable group and $M$ is
$(n)$-solvable. If $\phi : \pi_1(M)\longrightarrow\Gamma$ extends
over some $(n.5)$-solution $W$ for $M$, then $\rho(M,\phi) =0$.
\end{thm}

\begin{cor}
\label{cor:rho-slice} If $K$ is a slice knot and $\phi$ extends
over the exterior of a slice disk, then $\rho(M,\phi) = 0$ for any
PTFA group $\Gamma$ where $M$ is zero surgery on $K$ in $S^3$.
\end{cor}

As to calculating $\rho$-invariants, if $\Gamma = \bbz$ and $\phi$
is not trivial , then $\rho(M,\phi)$ is easily calculated as a
certain integral over $S^1$. (See \cite[Property 2.4] {COT2}.) In
particular, if $K$ has Arf invariant zero (i.e., $K$ is
$(0)$-solvable), then we can assign a real value $\rho(K)$ to $K$
that is ``canonically" induced from $\rho$-invariants as follows.
Let $M$ be zero surgery on $K$ in $S^3$. Choose a $(0)$-solution
$W$ of $M$ such that a canonical epimorphism $\phi :
\pi_1(M)\longrightarrow\bbz$ extends to $\psi : \pi_1(W)
\longrightarrow \bbz$ and $\pi_1(W) \cong \bbz$ as we did in the
proof of Proposition~\ref{prop:n-solvability}. Then we can
calculate $\rho(M,\phi)$ via $(W,\psi)$ and we define $\rho(K)$ to
be $\rho(M,\phi)$. These ``canonical" real numbers will play an
important role in our work. (See
Proposition~\ref{prop:n-solvability} and the paragraph preceding
Proposition~\ref{prop:1-solvability}.)

Now we investigate how $\rho$-invariants change under the grafting
construction described in Section~\ref{sec:n-solvable knots}.
Though there is no general additive property of $\rho$-invariants,
if the representations of the fundamental groups of the relevant
manifolds are matched up nicely under the grafting construction,
we can derive an additive property. In particular, to prove the
main theorem, we only need to look into $\rho$-invariants of
$K(J,\eta)$ where $K$ is a ribbon knot and $J$ has Arf invariant
zero.

Suppose $K$ is a ribbon knot and $J$ has Arf invariant zero. Let
$W$ be the exterior of a ribbon disk for $K$. $W_J$, $W'$, $M_J$,
$M'$, $\eta$ are defined as in
Proposition~\ref{prop:n-solvability}. Suppose we are given
homomorphisms $\phi : \pi_1(M)\longrightarrow\Gamma$ and $\phi_J :
\pi_1(M_J)\longrightarrow\Gamma$ such that $\phi([\eta]) =
\phi_J([\mu_J])$ where $\Gamma$ is a PTFA group. Then $\phi$ and
$\phi_J$ produce a unique homomorphism $\phi' :
\pi_1(M')\longrightarrow\Gamma$ (For this, observe that $M' =
(M\setminus (\eta\times D^2))\bigcup\limits_{S^1\times S^1}E_J$
where $M\setminus (\eta\times D^2)\subset M$ and $E_J\subset M_J$.
Use Van Kampen Theorem  noticing that for $\{*\}\times\partial
D^2\subset\eta\times D^2$, $\phi([\{*\}\times\partial D^2]) =
\phi_J([\ell_J]) = 0$). Then we have the following proposition due
to COT.

\begin{prop}
\label{prop:additivity of rho invariant} Suppose $\phi'$ extends
to $\psi' : \pi_1(W')\longrightarrow\Gamma$ . Then $\rho(M',\phi')
= \rho(J)$ if $\phi(\eta)\not=1$, and $\rho(M',\phi') = 0$ if
$\phi(\eta) = 1$.
\end{prop}
\begin{proof}
By \cite[Proposition 3.2]{COT2}, $\rho(M',\phi') = \rho(M,\phi) +
\rho(M_J, \phi_J)$. $\rho(M,\phi) = 0$ by
Corollary~\ref{cor:rho-slice}. Since $\eta$ generates
$\pi_1(W_J)\cong\bbz$, $\rho(M_J,\phi_J) = \rho(J)$ if $\phi(\eta)
\not= 1$ and $\rho(M_J, \phi_J) = 0$ if $\phi(\eta) = 1$ by
Property 2.3 and 2.5 in \cite{COT2}.
\end{proof}

\section{Proof of main theorem}
\label{sec:main theorem} We begin this section by briefly
reviewing some very useful machinery for the proof of the main
theorem. This originates from \cite[Section 2, 3, and 4]{COT1}, so
for all the detailed arguments and more generalized facts, the
readers should consult \cite{COT1}.

\begin{definition}
\label{def:univ gr} The family of {\em rationally universal}
groups $\{\gu\}$ is defined inductively by $\Gamma^U_0 = \bbz$,
$\scrr^U_0=\bbq[t,t^{-1}]$ and for $n\geq 0$, setting
\[
S_n=\bbq[\gu,\gu]-\{0\}, \phantom{.......}
\scrr^U_n=(\bbq\gu)S^{-1}_n
\]
and
\[
\Gamma^U_{n+1}=\scrk_n/\scrr^U_n\rtimes\gu.
\]
Here $\scrk_n$ is the right ring of quotients of $\bbq[\gu]$.
\end{definition}
\noindent It is shown in \cite[Proposition 2.5]{COT1} that
$\bbq[\gu]$ is an Ore domain, i.e., $\bbq[\gu]$ has a right ring
of quotients. (Note that inductively $\gu$ is PTFA.)
The semi-direct product is defined via the left multiplication of
$\gu$ on $\scrk_n/\scrr^U_n$. One can show that $\gu$ is an
$(n)$-solvable  group for all $n\ge 0$. Observe that $\scrk_0 =
\bbq(t)$ and $\Gamma^U_1 = \bbq(t)/\bbq[t,t^{-1}]\rtimes\bbz$.

Suppose $M$ is a closed 3-manifold with $\beta_1(M) = 1$ and we
have a homomorphism $\phi_{0} :
\pi_1(M)\longrightarrow\Gamma^U_{0}$. Then we can define the
rational Alexander module $\scra_{0}(M)\equiv H_1(M;\scrr^U_{0})$
and the (non-singular) Blanchfield form $B\ell_0 :
\scra_0(M)\times\scra_0(M) \longrightarrow \scrk_0/\scrr^U_0$.
Then,
\[
\scra_{0}(M)\equiv H_1(M;\scrr^U_{0})\cong H^2(M;\scrr^U_{0})\cong
H^1(M;\scrk_{0}/\scrr^U_{0}).
\]
and there is a bijection $f : H^1(M;\scrk_{0}/\scrr^U_{0})
\longleftrightarrow
\mbox{Rep}^{\ast}_{\Gamma^U_{0}}(\pi_1(M),\Gamma^U_1)$
\newline\noindent ($\mbox{Rep}^{\ast}_{\Gamma^U_{0}}
(\pi_1(M),\Gamma^U_1)$ is defined to be the representations from
$\pi_1(M)$ to $\Gamma^U_1$ which agree with $\phi_{0}$ after
composing with the projection $\Gamma^U_1
\longrightarrow\Gamma^U_{0}$ modulo
$\scrk_{0}/\scrr^U_{0}$-conjugations.) So any choice
$x_{0}\in\scra_{0}(M)$ will (together with $\phi_{0}$) induce
$\phi_1 : \pi_1(M)\longrightarrow\Gamma^U_1$. We refer to this as
{\em the coefficient system corresponding to $x_{0}$ (and
$\phi_{0}$)}. One can think of (the image of) this element $x_{0}$
as an element of
$\Hom_{\scrr^U_{0}}(\scra_{0}(M),\scrk_{0}/\scrr^U_{0})$ under the
Kronecker map from $H^1(M;\scrk_{0}/\scrr^U_{0})$. This image is
called the {\em character induced by} $x_{0}$. Now we obtain some
very useful facts which are summarized in the following remark.

\begin{remark}
\label{rem:useful facts} (\cite[Theorem 3.5, 3.6, and 4.4]{COT1})
Suppose $M = \partial W$ is a compact 3-manifold with $\beta_1(M)
= 1$ and $\phi_0 : \pi_1(M) \longrightarrow \Gamma^U_0$ is given.
\begin{itemize}
\item [(i)] The isomorphism $H_1(M;\scrr^U_0)\cong H^1(M;\scrk_0/\scrr^U_0)$ with
$f$ gives a natural bijection $\tilde f : \scra_0(M)
\longleftrightarrow \mbox{Rep}^{\ast}_{\Gamma^U_0}
\left(\pi_1(M),\Gamma^U_{1}\right)$.
\item [(ii)] If $x\in\scra_0(M)$, then the character induced by $x$ is
given by $y\mapsto B\ell_0(x,y)$.
\item [(iii)] Assume that the non-trivial map
$\phi_{0}:\pi_1(M)\longrightarrow\Gamma^U_{0}$ extends to a map
$\psi_{0}:\pi_1(W)\longrightarrow\Gamma^U_{0}$ and that $\phi_1$
is a representative of a class in
$\mbox{Rep}^{\ast}_{\Gamma^U_{0}}\left(\pi_1(M),\Gamma^U_1\right)$
corresponding to $x\in H_1(M;\scrr^U_{0})$. Let
\[P_{0}\equiv\mbox{$\Ker$}\{j_{\ast}:H_1(M;\scrr^U_{0})\longrightarrow
H_1(W;\scrr^U_{0})\}.\] Then if $M$ is $(1)$-solvable via $W$,
then $\phi_1$ extends to $\pi_1(W)$ if and only if $x\in P_{0}$.
\item [(iv)] Suppose $M$ is $(1)$-solvable via $W$ and $\phi_0$
is a non-trivial coefficient system that extends to $\pi_1(W)$.
Then the Blanchfield form $B\ell_0$ is hyperbolic, and in fact the
kernel of $j_{\ast}:H_1(M;\scrr^U_0) \longrightarrow
H_1(W;\scrr^U_0)$ is self-annihilating. (i.e., $\ker j_{\ast} =
(\ker j_{\ast})^\bot$.)
\end{itemize}
\end{remark}

From Theorem~\ref{thm:facts of rho invariant} and
Remark~\ref{rem:useful facts}, to prove that a knot $K$ is not
$(1.5)$-solvable, basically we need to investigate the
representations of the fundamental group induced from all
self-annihilating submodules of $\scra_0(M)\equiv
H_1(M;\bbq[t,t^{-1}])$ where $M$ is zero surgery on $K$ in $S^3$.
But in case $\scra_0(M)$ has a unique proper submodule, we have
the following useful lemma.

\begin{lem}
\label{lem:uniqueness} Suppose $M$ is $(1)$-solvable and
$\scra_0(M)$ has a unique proper submodule $P$. If there exists
$p\in P$ such that $\rho(M,\phi)\not= 0$ for $\phi :
\pi_1(M)\longrightarrow \Gamma^U_1$ induced from $p$, then $M$ is
not $(1.5)$-solvable.
\end{lem}
\begin{proof}
Suppose $M$ is $(1.5)$-solvable via $W$. Let $\scra_0(W) \equiv
H_1(W;\bbq[t,t^{-1}])$. Since $W$ is also a $(1)$-solution of $M$,
the kernel of the inclusion-induced map $i_{\ast} : \scra_0(M)
\longrightarrow \scra_0(W)$ is self-annihilating with respect to
$B\ell_0$ by Remark \ref{rem:useful facts}~(iv). Since
$\scra_0(M)$ has a unique proper submodule $P$, $\Ker i_{\ast} =
P$. By Remark \ref{rem:useful facts}~(i) and (iii) ,
$\phi:\pi_1(M)\longrightarrow\Gamma_1^U$ induced from $p(\in P)$
extends to $\pi_1(W)$. Then since $W$ is assumed to be a
$(1.5)$-solution of $M$, by Theorem~\ref{thm:facts of rho
invariant}, $\rho(M,\phi) = 0$. This leads us to a contradiction.
\end{proof}

\begin{figure}[htb]
\includegraphics[scale=0.9]{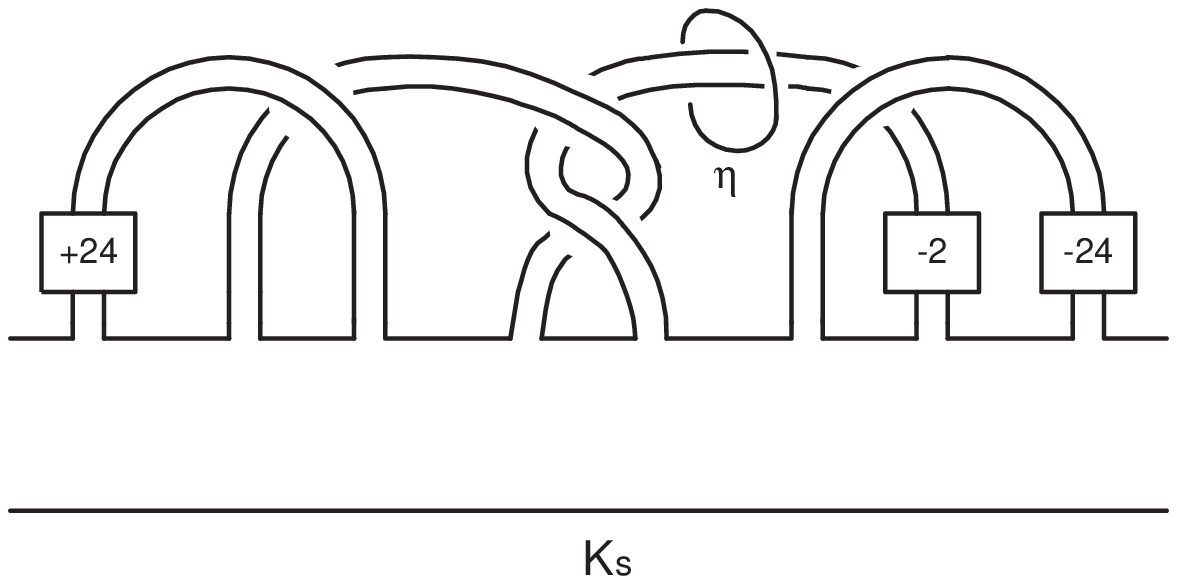}
\caption{}\label{fig:eta}
\end{figure}

Through this section, $K_s$ denotes our seed ribbon knot which was
constructed in Section~\ref{sec:seed knot} and $\eta$ is the
designated circle in the complement of the Seifert surface of
$K_s$ in $S^3$ as in Figure~\ref{fig:eta}. Notice that $\eta$ is a
dual of $w_9$, so it represents the homology class which generates
the rational Alexander module of $K_s$.

By \cite[Proposition 2.6]{COT2}, there are infinitely many Arf
invariant zero knots $J_i(i\in\bbn)$ such that
$\{\rho(J_i)\}_{i\in\bbn}$ is linearly independent over integers.
In particular, $\rho(J_i)\not=0$. Let $K_i \equiv K_s(J_i,\eta)$
be the family of knots resulting from the grafting construction as
described in Section~\ref{sec:n-solvable knots}. In the following
propositions, we exploit important properties of $K_i$.

\begin{prop}
\label{prop:1-solvability} $K_i$ $(i\in\bbn)$ are $(1)$-solvable
but not $(1.5)$-solvable.
\end{prop}
\begin{proof}
$\eta$ lifts to a closed circle in the infinite cyclic cover of
$S^3\setminus K_s$, hence $\eta\in\pi_1(M)^{(1)}$ where $M$ is
zero surgery on $K_s$ in $S^3$. Now it is clear from
Proposition~\ref{prop:n-solvability} that $K_i$ are
$(1)$-solvable.

We need to show that $K_i$ are not $(1.5)$-solvable. Fix $i$. Let
$W'$ denote the $(1)$-solution for $K_i$ formed as in the proof of
Proposition~\ref{prop:n-solvability} and let $M'$ denote zero
surgery on $K_i$ in $S^3$. Recall that $\Gamma^U_0 = \bbz$. Let
$\pi_1(M')\longrightarrow \Gamma^U_0$ be the canonical epimorphism
which extends uniquely to an epimorphism $\pi_1(W')$. Looking into
the grafting construction more closely, one can see that $K_i$ has
the same Seifert form as that of $K_s$, so the rational Alexander
module $\scra_0(M')$ is isomorphic to
$\bbq[t,t^{-1}]/(\Phi_{30}(t))^2$. Let $\scra_0(W') \equiv
H_1(W';\bbq[t,t^{-1}])$. By Remark \ref{rem:useful facts}~(iv),
since $W'$ is a $(1)$-solution for $M'$, the kernel of the
inclusion-induced map $i_* : \scra_0(M')\longrightarrow
\scra_0(W')$ is self-annihilating with respect to the
(non-singular) Blanchfield form $B\ell_0$. Since $\scra_0(M')$ has
a unique proper submodule, say $P_0$, which is generated by
$\Phi_{30}(t)$, $\Ker i_{\ast} = P_0$. Choose a non-zero $p_0\in
P_0$ such that $B\ell_0(\eta,p_0)\not= 0$. Such a $p_0$ exists
since $\eta$ generates $\scra_0(M')$ and $B\ell_0$ is non-singular
for which $P$ is self-annihilating. Then $p_0$ induces $\phi :
\pi_1(M') \longrightarrow \Gamma_1^U$ by Remark \ref{rem:useful
facts}~(i). By Remark \ref{rem:useful facts}~(iii), $\phi$ extends
to $\psi : \pi_1(W') \longrightarrow \Gamma_1^U$. Now we compute
$\rho(M',\phi)$ using $(W',\psi)$. Since $B\ell_0(\eta,p_0)\not=
0$, $\phi(\eta)\not= 1$ by Remark \ref{rem:useful facts}~(ii). By
Proposition~\ref{prop:additivity of rho invariant}, $\rho(M',\phi)
= \rho(J_i)$, which is nonzero by our choice of $J_i$. By
Lemma~\ref{lem:uniqueness}, $K_i$ is not $(1.5)$-solvable.
\end{proof}

\begin{prop}
\label{prop:vanishing CG}
 $K_i$ have vanishing Casson-Gordon invariants.
\end{prop}
\begin{proof}
Because $K_i$ and $K_s$ have the same Seifert form, they have the
same Alexander polynomial which is $(\Phi_{30}(t))^2$. By
Theorem~\ref{thm:homology sphere}, any prime power branched cyclic
cover of $K_i$ is a homology sphere. Hence all Casson-Gordon
invariants vanish on $K_i$ by \cite[Corollary B2]{Lit}.
\end{proof}

Now we are ready to prove the main theorem.

\begin{proof}[Proof of Theorem~\ref{thm:main theorem}]
First, we show that no non-trivial linear combination of $K_i,
i\in \bbn$ is $(1.5)$-solvable. We follow \cite{COT2}. Refer to
the proof of \cite[Theorem 4.1]{COT2}. The only crucial difference
in this proof is that we deal with $(1.5)$-solvability instead of
$(2.5)$-solvability, so we use the second order invariants instead
of the third order invariants. Since the proof will follow almost
the same course of COT's proof for \cite[Theorem 4.1]{COT2}, some
details will be omitted. For convenience, we follow the notations
used in \cite{COT2}.

Suppose that a non-trivial linear combination $\#^m_{i=1}n'_iK_i,
n'_i\not= 0$, is $(1.5)$-solvable. We may assume all $n'_i > 0$ by
replacing $K_i$ by $-K_i$ if $n'_i < 0$ and  $n'_1 > 1$ if $m =
1$. Let $M_i$ denote $M_{K_i}$, and note that $-M_i = M_{-K_i}$.
Let $M_0$ denote $0$-surgery on $\#^m_{i=1}n'_iK_i$. Let $W_0$ be
a $(1.5)$-solution of $M_0$. Let $W_i$ $(i
> 0)$ denote the specific $(1)$-solution for $M_i$ constructed as
in Proposition~\ref{prop:n-solvability} with the exterior of a
ribbon disk for $K$. Let $n_1 = n'_1 - 1$ and $n_i = n'_i$ if $i >
1$. Let $W$ be the union of $W_0$, $C$ ($C$ is defined in the next
paragraph), and all the copies of $W_i$ $(i=1, 2, \dots, m)$ where
there are $n_i$ copies of $W_i$ below $C$. Refer to
Figure~\ref{fig:cobordism} below. Later we will show that $W$ is a
$(1)$-solution of $M_1$.

\begin{figure}[htb]
\includegraphics[scale=1]{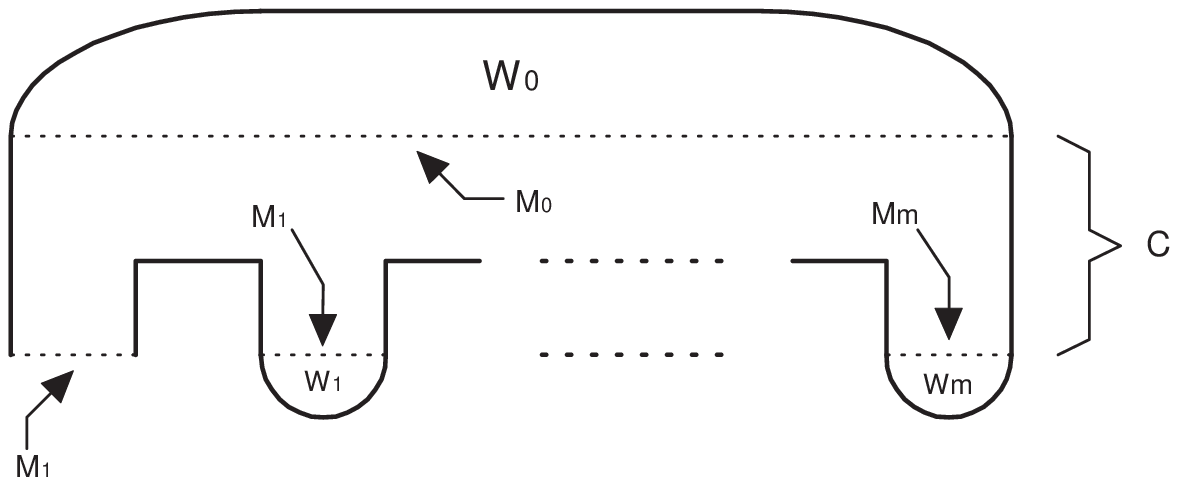}
\caption{}\label{fig:cobordism}
\end{figure}

The 4-manifold $C$ is a standard cobordism between $0$-surgery on
$\#^m_{i=1}n'_iK_i$ and the disjoint union of $0$-surgeries on the
summands of $\#^m_{i=1}n'_iK_i$. Briefly, start with a collar on
the disjoint union of $0$-surgeries, and add $1$-handles to get a
connected 4-manifold whose upper boundary is given by surgery on
the link consisting of the split union of $K_i$'s, each with
$0$-framing. Next add $0$-framed $2$-handles to get a $0$-surgery
on a connected sum of $K_i$'s on the upper boundary. See
\cite[Theorem 4.1]{COT2} for more details, and note that $C$ has a
handlebody decomposition, relative to $\coprod^m_{i=1}n'_iM_i$,
consisting of $(\sum^m_{i=1}|n_i|)$ $1$-handles and the same
number of $2$-handles. Moreover, $H_1(C;\bbz) \cong \bbz$ and the
inclusion from any of its boundary components induces an
isomorphism on $H_1$. One can also see that $H_2(C)\cong
H_2(\coprod n'_iM_i)$.

We prove that $W$ is a $(1)$-solution of $M_1$. Since the
inclusion-induced homomorphisms $H_1(M_i)\longrightarrow H_1(W_i)$
are isomorphisms for $i\ge 0$, the inclusion-induced
$H_1(M_1)\longrightarrow H_1(W)$ is also an isomorphism.  For $i
\ge 0$, $H_2(M_i)\longrightarrow H_2(W_i)$ is the zero map since
the boundary map $H_3(W_i,M_i)\longrightarrow H_2(M_i)$, the dual
map of the inclusion induced $H^1(W_i)\longrightarrow H^1(M_i)$,
is an isomorphism. Using this and Mayer-Vietoris sequence  we can
prove that $H_2(W) \cong H_2(W_0)\bigoplus^m_{i=1}n_iH_2(W_i)$.
Now if one looks carefully at $(1)$-Lagrangians and their duals
for $W_0$ and the $W_i$'s, one can see that they form
$(1)$-Lagrangian and its dual for $W$. So $W$ is a $(1)$-solution
for $M_1$.

We repeat the argument in Proposition~\ref{prop:1-solvability}.
Let $\pi_1(M_1)\longrightarrow \Gamma^U_0$ be the canonical
epimorphism which extends uniquely to an epimorphism $\pi_1(W)$.
Recall that the rational Alexander module $\scra_0(M_1) =
H_1(M_1;\bbq[t,t^{-1}])$ is isomorphic to
$\bbq[t,t^{-1}]/(\Phi_{30}(t))^2$. Let $\scra_0(W) =
H_1(W;\bbq[t,t^{-1}])$. By Remark \ref{rem:useful facts} (iv),
since $W$ is a $(1)$-solution for $M_1$, the kernel of the
inclusion-induced map $j_* : \scra_0(M_1) \longrightarrow
\scra_0(W)$ is self-annihilating with respect to the Blanchfield
form $B\ell_0$. Since $\scra_0(M_1)$ has a unique proper
submodule, say $P_0$, the latter is this kernel. Choose a non-zero
$p_0\in P_0$, inducing $\phi_1 : \pi_1(M_1) \longrightarrow
\Gamma^U_1$ by Remark \ref{rem:useful facts} (i). By Remark
\ref{rem:useful facts} (iii), $\phi_1$ extends to $\psi_1 :
\pi_1(W) \longrightarrow \Gamma^U_1$. Therefore $\rho(M_1,\phi_1)$
can be computed using $(W,\psi_1)$.

We compute $\rho(M_1,\phi_1)$ using $(W,\psi_1)$. Let
$\phi_{(i,j)}$ denote the restriction of $\psi_1$ to the
$j^{\text{th}}$ copy of $\pi_1(M_i), 1\le j\le n_i$. Let $\phi_0$
denote the restriction of $\psi_1$ to $\pi_1(M_0)$. Let $\scrk_1$
denote the classical right ring of quotients of $\bbz\Gamma_1^U$.
$H_*(M_i;\scrk_1) = 0$ for $i\ge 0$ (\cite[Propositions 2.9 and
2.11]{COT1}), so a Mayer-Vietoris sequence shows that
$H_2(W;\scrk_1) \cong H_2(W_0;\scrk_1)\oplus H_2(C;\scrk_1)\oplus
H_2(W_1;\scrk_1)\oplus\cdots\oplus H_2(W_m;\scrk_1)$ where $W_i$
occurs $n_i$ times. Here the coefficient systems on $W_i$ $(i\ge
0)$ and $C$ are induced by inclusions into $W$. By \cite[Lemma
4.2]{COT2} $H_2(C;\scrk_1) = 0$. And the intersection form on
$H_2(W;\scrk_1)$ splits along the direct sum.  From \cite[Section
5]{COT1}, $\sigma^{(2)}_{\Gamma^U_1}$ can be viewed as a
homomorphism from the Witt group of non-singular hermitian forms
on finitely generated $\scrk_{1}$ modules, so we have

\[
 \rho(M_1,\phi_1) = \rho(M_0,\phi_0) +
 \sum^m_{i=1}\sum^{n_i}_{j=1}\rho(-M_i,\phi_{(i,j)})
\]

Here $\rho(M_0,\phi_0) = 0$ by Theorem~\ref{thm:facts of rho
invariant} because $\phi_0$ extends to $(1.5)$-solution $W_0$. By
Proposition~\ref{prop:additivity of rho invariant},
$\rho(-M_i,\phi_{(i,j)}) = -\rho(M_i,\phi_{(i,j)}) = -\rho(J_i)$
or 0. So we deduce that

\[
 \rho(M_1,\phi_1) + \sum^m_{i=1}c_i\rho(J_i) = 0
\]

\noindent for some non-negative constants $c_i$'s.

Now as in Proposition~\ref{prop:1-solvability}, pick $p_0 \in P_0$
such that $\phi_1(\eta) \not= 1$. Note that $P_0$ is equal to the
kernel of the inclusion-induced map $i_* :
\scra_0(M_1)\longrightarrow \scra_0(W_1)$. Then by
Proposition~\ref{prop:additivity of rho invariant},
$\rho(M_1,\phi_1) = \rho(J_1)$, so we have

\[
 \rho(J_1) + \sum^m_{i=1}c_i\rho(J_i) = 0
\]
\noindent which contradicts that $\{\rho(J_i)\}_{i\in\bbn}$ is
linearly independent over integers.

Now it remains to show that Casson-Gordon invariants vanish on the
subgroup generated by $K_i$. Recall that every prime power
branched cyclic cover of $K_i$ is a homology 3-sphere. One can
show that the connected sum of two homology 3-spheres is a
homology 3-sphere. Since a finite branched cyclic cover of the
connected sum of two knots over $S^3$ is homeomorphic with the
connected sum of the finite branched cyclic covers of the knots,
the assertion follows from \cite[Corollary B2]{Lit}.
\end{proof}

\end{document}